# SQUARE ROOT PENALTY: ADAPTATION TO THE MARGIN IN CLASSIFICATION AND IN EDGE ESTIMATION

By A. B. Tsybakov and S. A. van de Geer

*Université Paris VI and University of Leiden*

We consider the problem of adaptation to the margin in binary classification. We suggest a penalized empirical risk minimization classifier that adaptively attains, up to a logarithmic factor, fast optimal rates of convergence for the excess risk, that is, rates that can be faster than $n^{-1/2}$, where $n$ is the sample size. We show that our method also gives adaptive estimators for the problem of edge estimation.

**1. Introduction.** Consider observations $(X_1, Y_1), \ldots, (X_n, Y_n)$, where $Y_i$ is a bounded response random variable and $X_i \in \mathcal{X}$ is the corresponding instance. We regard $\{(X_i, Y_i)\}_{i=1}^n$ as i.i.d. copies of a population version $(X, Y)$. The goal is to predict the response $Y$ given the value of the instance $X$. We consider two statistical problems: binary classification and boundary estimation in binary images (edge estimation). In the classification setup $Y_i \in \{0, 1\}$ is a label (e.g., {ill, healthy}, {white, black}, etc.), while in edge estimation $Y_i$ can be either a label or a general bounded random variable. Most of the paper will be concerned with the model of binary classification. The results for edge estimation are quite analogous and they will be stated as corollaries in Section 6.

Any subset $G$ of the instance space $\mathcal{X}$ may be identified with its indicator function $\mathbb{1}_G$, that is, with a *classification rule* or *classifier* $G$ which predicts $Y = 1$ iff $X \in G$. The prediction error $R(G)$ of the classifier $G$ is the probability that it predicts the wrong label, that is,

$$(1.1) \qquad R(G) = \mathbf{E}([Y - \mathbb{1}_G(X)]^2).$$

Let $\eta(X) = P(Y = 1 | X)$ be the regression of $Y$ on $X$. The Bayes rule is the classifier

$$(1.2) \qquad G^* = \{x \in \mathcal{X} : \eta(x) > 1/2\}.$$

This rule is optimal in the sense that it minimizes the prediction error over all $G \subset \mathcal{X}$ [see, e.g., Devroye, Györfi and Lugosi (1996)]. The regression $\eta$ is generally unknown. We consider the construction of an estimator $\hat{G}_n \subset \mathcal{X}$ of the Bayes rule $G^*$ without directly estimating $\eta$.

The performance of a classifier $\hat{G}_n$ is measured by its excess risk $\mathbf{E}(R(\hat{G}_n)) - R(G^*)$. It is well known that for various classifiers the excess risk converges to 0 as $n \to \infty$ at the rate $n^{-1/2}$ or slower [see Devroye, Györfi and Lugosi (1996) and Vapnik (1998), where one can find further references]. Moreover, under conditions on the identifiability of the minimum of the risk $R(\cdot)$ called *margin conditions*, some classifiers can attain *fast rates*, that is, rates that are faster than $n^{-1/2}$. The existence of such fast rates in classification problems has been established by Mammen and Tsybakov (1999). They showed that optimal rates of convergence of the excess risk to 0 depend on two parameters: complexity of the class of candidate sets $G$ (parameter $\rho$) and the *margin parameter* $\kappa$ which characterizes the extent of identifiability. Their construction was nonadaptive supposing that $\rho$ and $\kappa$ were known. Tsybakov (2004) suggested an adaptive classifier that attains the fast optimal rates, up to a logarithmic factor, without prior knowledge of the parameters $\rho$ and $\kappa$, thus solving the so-called *adaptation to the margin* problem. The classification rule suggested by Tsybakov (2004) is based on multiple pretesting aggregation of empirical risk minimizers over a collection of classes of candidate sets $G$. This procedure differs significantly from penalized empirical risk classifiers that are widely used in modern practice of classification [cf. Schölkopf and Smola (2002)]. Subsequently there has been a discussion in the literature of whether penalized classifiers can adaptively attain fast optimal rates. In particular, Koltchinskii and Panchenko (2002) and Audibert (2004) proposed convex combinations of classifiers, and Koltchinskii (2001) and Lugosi and Wegkamp (2004) suggested data-dependent penalties. The resulting adaptive classifiers converge with rates that can be faster than $n^{-1/2}$ but that are different from the optimal rates in a minimax sense considered in Tsybakov (2004).

This paper answers affirmatively to the above question: penalized classifiers can adaptively attain fast optimal rates. Moreover, the penalty allowing one to achieve this effect is not data-dependent or randomized. It is very simple and essentially arises from a sparsity argument similar to the one used in the wavelet thresholding context. Interestingly, the penalty is not of the $\ell_1$-type as for soft thresholding and not of the $\ell_0$-type as for hard thresholding, but rather of an intermediate, block-wise $\ell_{1/2}$ or "square root" type. Inspection of the proof shows that the effect is very pointed, that is, the proof heavily relies on our particular choice of the penalty.

The classifier $\hat{G}_n$ that we study is constructed as follows. Let

$$R_n(G) = \frac{1}{n} \sum_{i=1}^{n} (Y_i - \mathbb{1}_G(X_i))^2 \tag{1.3}$$

be the empirical risk of a classifier $G \subset \mathcal{X}$. Note that $R_n(G)$ is the proportion of observations misclassified by $G$ and that its expectation $R(G) = \mathbf{E}(R_n(G))$ is the prediction error. Assume that $X = (S,T) \in \mathcal{X} = [0,1]^{d+1}$, with $S \in [0,1]^d$ ($d \leq \log n$), and $T \in [0,1]$. A boundary fragment is a subset $G$ of $\mathcal{X}$ of the form

$$G = \{(s,t) \in \mathcal{X} : f(s) \geq t\} \tag{1.4}$$

where $f$ is a function from $[0,1]^d$ to $[0,1]$ called the *edge function*. We let $\hat{G}_n$ be a minimizer of the penalized empirical risk

$$R_n(G) + \operatorname{Pen}(G) \tag{1.5}$$

over a large set of boundary fragments $G$. Here $\operatorname{Pen}(G)$ is a penalty on the roughness of the boundary. The purpose of the penalty is to avoid overfitting. We will show that a weighted *square root penalty* [see (2.2) and (2.3)] results in a classifier with the adaptive properties as discussed above.

A refinement as compared to Tsybakov (2004) is that we do not only consider adaptation in a minimax sense but also adaptation to the oracle. We obtain asymptotically exact oracle inequalities and then get minimax adaptation as a consequence. We work under somewhat different assumptions than in Tsybakov (2004). They are slightly more restrictive as concerns the model. For example, we consider only boundary fragments as candidates for $G$. The class of boundary fragments is possibly a genuine restriction, although some generalizations to other classes of sets are clearly feasible. On the other hand, our assumptions allow us to adapt to more general smoothness (complexity) properties of $G$. For example, Vapnik–Chervonenkis classes of sets $G$ (corresponding approximately to $\rho = 0$, see Section 5) or the classes of sets with very nonsmooth boundaries (corresponding to $\rho \geq 1$) are covered by our approach.

As a corollary of the results, we obtain an adaptive estimator in the problem of edge estimation considered by Korostelev and Tsybakov (1993). The statistical model in that problem is similar to the one described above. However, it treats the situation characteristic for image analysis where the $X_i$'s are uniformly distributed on $\mathcal{X}$, and the error criterion is not the excess risk but rather the risk $\mathbf{E}(\mu_{d+1}(\hat{G}_n \triangle G^*))$, where $\triangle$ is the symbol of symmetric difference between sets and $\mu_{d+1}$ denotes the Lebesgue measure on $[0,1]^{d+1}$.

The paper is organized as follows. In Section 2 we define our adaptive classifier. In Section 3 we introduce some notation and assumptions. Section 4 presents the main oracle inequality. In Section 5 we apply this inequality to get minimax adaptation results. Section 6 discusses the consequences for edge estimation. Proofs are given in Section 7.

**2. Definition of the adaptive classifier.** Let $\{\psi_k : k = 1, \ldots, n\}$ be an orthonormal system in $L_2([0,1]^d, \mu_d)$ where $\mu_d$ is the Lebesgue measure on $[0,1]^d$. For $\alpha \in \mathbf{R}^n$ define

$$f_\alpha(s) = \sum_{k=1}^n \alpha_k \psi_k(s), \qquad s \in [0,1]^d. \tag{2.1}$$

Introduce a double indexing for the system $\{\psi_k\}$, namely

$$\{\psi_k : k = 1, \ldots, n\} = \{\psi_{j,l} : j \in I_l, l = 1, \ldots, L\}$$

where $I_l$, $l = 1, \ldots, L$, are disjoint subsets of $\{1, \ldots, n\}$ such that

$$\sum_{l=1}^L |I_l| = n.$$

Here $|A|$ denotes the cardinality of the set $A$. One may think of $\{\psi_{j,l}\}$ as of a wavelet-type system with the index $l$ corresponding to a resolution level. A vector $\alpha \in \mathbf{R}^n$ can be written with this double indexing as $\alpha = (\alpha_{j,l})$.

For a linear classification rule defined by the set $G_\alpha = \{(s,t) \in \mathcal{X} : f_\alpha(s) \geq t\}$, consider the penalty

$$\text{Pen}(G_\alpha) = \lambda_n \sqrt{I(\alpha)}, \tag{2.2}$$

where $I(\cdot)$ is a nonsparsity measure of the form

$$I(\alpha) = \left(\sum_{l=1}^L w_l^{1/2} \sqrt{\sum_{j \in I_l} |\alpha_{j,l}|}\right)^2, \tag{2.3}$$

for certain weights $(w_l)$. In what follows we take the weights as

$$w_l = 2^{dl/2}, \qquad l = 1, \ldots, L, \tag{2.4}$$

and we prove our results for wavelet-type bases (cf. Assumption B below). An extension to other bases $\{\psi_k\}$ is possible where the block sizes $|I_l|$ should be chosen in an appropriate way [e.g., as in Cavalier and Tsybakov (2001)]. The weights $w_l$ should moreover be defined as a function of $|I_l|$. We do not pursue this issue here because it requires different techniques. Thus, in this paper we consider penalties based on

$$I(\alpha) = \left(\sum_{l=1}^L 2^{dl/4} \sqrt{\sum_{j \in I_l} |\alpha_{j,l}|}\right)^2.$$

One may think of $\{\alpha_{j,l}\}$ as the coefficients of the expansion of a function in the Besov space $B_{\sigma,p,q}([0,1]^d)$, with $p = 1$, $q = 1/2$ and smoothness $\sigma = (d+1)d/2$ [so that the effective smoothness is $s = \sigma/d = (d+1)/2$; see, for example, DeVore and Lorentz (1993).

We propose the estimator $\hat{G}_n = G_{\hat{\alpha}_n}$ where

(2.5) $$\hat{\alpha}_n = \arg\min_{\alpha \in \mathbf{R}^n} \{R_n(G_\alpha) + \lambda_n \sqrt{I(\alpha)}\}.$$

Here $\lambda_n > 0$ is a regularization parameter that will be specified in Theorem 1. We refer to $\lambda_n \sqrt{I(\alpha)}$ as a (block-wise) $\ell_{1/2}$ or *square root penalty*.

One may compare (2.5) to a wavelet thresholding estimator for regression. The difference here is that because our problem is nonlinear, we cannot express the solution $\hat{\alpha}_n$ in a levelwise form, and we need to treat all the coefficients $\alpha_{j,l}$ globally.

**3. Notation and assumptions.** Let $G \triangle G'$ be the symmetric difference between two sets $G$ and $G'$, and let $Q$ denote the distribution of $X$. For a Borel function $f : [0,1]^d \to [0,1]$, we let

(3.1) $$\|f\|_1 = \int |f(s)| \, d\mu_d(s)$$

be its $L_1$-norm. (Recall that $\mu_d$ denotes the Lebesgue measure on $[0,1]^d$.) Note that

(3.2) $$\mu_{d+1}(G_\alpha \triangle G_{\alpha'}) = \|f_\alpha - f_{\alpha'}\|_1 = \|f_{\alpha-\alpha'}\|_1.$$

ASSUMPTION A. For some (unknown) $\kappa \geq 1$ and $\sigma_0 > 0$ and for all $\alpha \in \mathbf{R}^n$ we have

(3.3) $$R(G_\alpha) - R(G^*) \geq \frac{1}{\sigma_0} Q^\kappa(G_\alpha \triangle G^*).$$

Assumption A is a condition on sharpness of identifiability for the minimum of the risk. We will call it the *margin condition*. We refer to Tsybakov (2004) for a discussion of this condition. In particular, it is related to the behavior of the probability $Q(|\eta(X) - 1/2| \leq t)$ for small $t$. The case $\kappa = 1$ corresponds to a jump of $\eta$ at the boundary of $G^*$, and this is the most favorable case for estimation, while $\kappa \to \infty$ corresponds to a "plateau" around the boundary, and this is the least favorable case. For more discussion of the margin condition in relation to convex aggregation of classifiers, such as boosting, see Bartlett, Jordan and McAuliffe (2003) and Blanchard, Lugosi and Vayatis (2003).

We will also require the following condition on the basis.

ASSUMPTION B. The system of functions $\{\psi_{j,l}, j \in I_l, l = 1, \ldots, L\}$ is orthonormal in $L_2([0,1]^d, \mu_d)$ and satisfies, for some constant $c_\psi \geq 1$,

(3.4) $$\|\psi_{j,l}\|_1 \leq c_\psi 2^{-dl/2}, \qquad l = 1, \ldots, L,$$

$$(3.5) \qquad \sup_{s \in [0,1]^d} \sum_{j \in I_l} |\psi_{j,l}(s)| \leq c_\psi 2^{dl/2}, \qquad l = 1, \ldots, L,$$

$$(3.6) \qquad 2^{dl}/c_\psi \leq |I_l| \leq c_\psi 2^{dl}$$

and

$$(3.7) \qquad L \leq c_\psi \frac{\log n}{d}.$$

Assumption B makes it possible to relate $\|f_\alpha\|_1$ to $I(\alpha)$ in a suitable way (cf. Lemmas 1 and 2). Note that Assumption B is quite standard. It is satisfied, for instance, for usual bases of compactly supported wavelets [cf. Härdle, Kerkyacharian, Picard and Tsybakov (1998), Chapter 7].

Note also that (3.7) follows from (3.6) with a different constant. To simplify the exposition and calculations, we take the same constant $c_\psi$ in all the conditions (3.4)–(3.7) and suppose that this constant is not smaller than 1.

REMARK 1. It will be clear from the proofs of Lemmas 1 and 2 that Assumption B can be relaxed. Namely, the orthonormality of $\{\psi_{j,l}\}$ and (3.4) can be replaced by the conditions

$$\sum_{j \in I_l} |\alpha_{j,l}| 2^{-dl/2} \leq c_\psi \|f_\alpha\|_1, \qquad l = 1, \ldots, L,$$

$$\|f_\alpha\|_1 / c_\psi \leq \sum_{l=1}^{L} \sum_{j \in I_l} |\alpha_{j,l}| 2^{-dl/2} \qquad \forall \alpha \in \mathbf{R}^n.$$

Finally we introduce an assumption which will allow us to interchange Lebesgue measure and $Q$.

ASSUMPTION C. The distribution $Q$ of $X$ admits a density $q(\cdot)$ with respect to Lebesgue measure in $[0,1]^{d+1}$, and for some constant $1 \leq q_0 < \infty$ one has $1/q_0 \leq q(x) \leq q_0$ for all $x \in [0,1]^{d+1}$.

**4. An oracle inequality.** For $\alpha \in \mathbf{R}^n$ let

$$(4.1) \qquad m(\alpha) = \min\{m : \alpha_{j,l} = 0 \text{ for all } j \in I_l \text{ with } l > m\}$$

and

$$(4.2) \qquad N(\alpha) = N_{m(\alpha)},$$

with

$$(4.3) \qquad N_m = \sum_{l=1}^{m} |I_l|, \qquad m = 1, 2, \ldots, L.$$

Assume that there exists $\alpha^{\text{oracle}} \in \mathbf{R}^n$ such that

$$
\begin{aligned}
R(G_{\alpha^{\text{oracle}}}) &- R(G^*) + V_n(N(\alpha^{\text{oracle}})) \\
&= \min_{\alpha \in \mathbf{R}^n} \{R(G_\alpha) - R(G^*) + V_n(N(\alpha))\},
\end{aligned}
\tag{4.4}
$$

where

$$
V_n(N) = 4c_\kappa (4c_d q_0 c_\psi^2 \sigma_0^{1/\kappa} \lambda_n^2 N)^{\kappa/(2\kappa-1)}
\tag{4.5}
$$

and where $c_\kappa = (2\kappa - 1)/(2\kappa)\kappa^{-1/(2\kappa-1)}$ and $c_d = 2(2^d - 1)/(2^{d/2} - 1)^2$. Note that $V_n(N(\alpha))$ depends on the regularization parameter $\lambda_n$, which we shall take of order $\sqrt{\log^4 n/n}$ [see (4.6) in Theorem 1 below]. Then $\alpha^{\text{oracle}}$ can be interpreted as an oracle attaining nearly ideal performance. In fact, the term $R(G_\alpha) - R(G^*)$ in (4.4) may be viewed as an approximation error, while $V_n(N(\alpha))$ is related to the stochastic error, as will be clear from the proofs. In other words, nearly ideal performance is attained by the value $\alpha^{\text{oracle}}$ that trades off generalized bias and variance.

THEOREM 1. *Suppose that Assumptions A–C are met. Then there exists a universal constant $C$ such that for*

$$
\lambda_n = C \sqrt{\frac{q_0 c_\psi^2 \log^4 n}{nd}}
\tag{4.6}
$$

*and for any $\delta \in (0,1]$ and $n \geq 8 q_0 c_\psi^2$, we have*

$$
\begin{aligned}
\mathbf{P}\bigg( R(\hat{G}_n) - R(G^*) &> (1+\delta)^2 \inf_{\alpha \in \mathbf{R}^n} \{R(G_\alpha) - R(G^*) \\
&\qquad + \delta^{-1/(2\kappa-1)} V_n(N(\alpha))\} + 2\lambda_n \sqrt{\frac{\log^4 n}{n}} \bigg) \\
&\leq C \exp\bigg[ -\frac{c_\psi \log^4 n}{C^2 d} \bigg].
\end{aligned}
\tag{4.7}
$$

Theorem 1 shows that, up to a constant factor and a small remainder term $2\lambda_n \sqrt{\log^4 n/n}$, the estimator $\hat{G}_n$ mimics the behavior of the oracle. If $\delta$ is chosen small enough or converging to 0, for example, $\delta = 1/\log n$, the factor preceding the infimum in the oracle inequality (4.7) approaches 1.

The regularization parameter $\lambda_n \asymp \sqrt{n^{-1} \log^4 n}$ appearing in Theorem 1 is larger than the choice $\sqrt{n^{-1} \log n}$ used for wavelet thresholding in regression or density estimation. The value of $\lambda_n$ is imposed by an inequality for the empirical process that controls the stochastic error. Lemma 4 presents such an inequality, and the additional $\log n$ factors are due to the result given there.

As a consequence of (4.7) and of the fact that $0 \leq R(G) \leq 1$ for all $G$, we get the following inequality on the excess risk:

$$\mathbf{E}(R(\hat{G}_n)) - R(G^*) \leq (1+\delta)^2 \inf_{\alpha \in \mathbf{R}^n} \{R(G_\alpha) - R(G^*) + \delta^{-1/(2\kappa-1)} V_n(N(\alpha))\}$$

$$+ 2\lambda_n \sqrt{\frac{\log^4 n}{n}} + C \exp\left[-\frac{c_\psi \log^4 n}{C^2 d}\right].$$

This inequality bounds the excess risk by the oracle risk of a linear classification rule $G_\alpha$ for *any* form of Bayes rule $G^*$. We emphasize that $G^*$ is not necessarily a boundary fragment, and $R(G_\alpha) - R(G^*)$ is not necessarily small. The results of this section are thus of the learning theory type [cf. Devroye, Györfi and Lugosi (1996) and Vapnik (1998)]. In the next section we will show that if $G^*$ is a boundary fragment satisfying some regularity conditions, the excess risk converges to zero at a fast rate.

**5. Minimax adaptation.** Here we will consider a minimax problem and we will show how the oracle inequality of Section 4 can be used to prove that our classifier adaptively attains fast optimal rates under smoothness assumptions on the edge function.

Since in a minimax setup results should hold uniformly in the underlying distribution, we first introduce some notation to express the dependence of the margin behavior on the distribution of $(X, Y)$. Let us keep $d$ and also $Q$ fixed. Then the joint distribution of $(X, Y)$ is determined by the conditional probability $\eta(x)$ of the event $Y = 1$ given that $X = x$. Let $\mathcal{H}$ be the class of all Borel functions $\eta$ on $\mathcal{X}$ satisfying $0 \leq \eta \leq 1$. For a given $\eta \in \mathcal{H}$, let $dP_\eta(x, y)$ be the probability measure

$$dP_\eta(x, y) = (y\eta(x) + (1-y)(1-\eta(x))) \, dQ(x), \qquad (x, y) \in \mathcal{X} \times \{0, 1\}.$$

Let $G_\eta^*$ be Bayes rule when $(X, Y)$ has distribution $P_\eta$. Finally, let $\mathbf{E}_\eta$ denote expectation w.r.t. the distribution of $\{(X_i, Y_i)\}_{i=1}^n$ under $P_\eta$. Now fix the numbers $\sigma_0 > 0$ and $\kappa \geq 1$ and define the collection of functions

$$\mathcal{H}_\kappa = \Big\{\eta \in \mathcal{H} : G_\eta^* = \{(s, t) \in \mathcal{X} : f_\eta^*(s) \geq t\},$$

(5.1) $$\frac{1}{\sigma_0} Q^\kappa(G_\alpha \triangle G_\eta^*) \leq R(G_\alpha) - R(G_\eta^*)$$

$$\leq \sigma_0 q_0^\kappa \|f_\alpha - f_\eta^*\|_\infty^\kappa \text{ for all } \alpha \in \mathbf{R}^n\Big\},$$

where $\|\cdot\|_\infty$ denotes the $L_\infty$-norm on $[0, 1]^d$ endowed with Lebesgue measure, and $R(\cdot)$ depends on $\eta$ but in the notation we omit this dependence for brevity. Note that we assume a lower as well as an upper bound for the excess risk in definition (5.1), and in view of Assumption C and (3.2), $Q^\kappa(G_\alpha \triangle G_\eta^*) \leq q_0^\kappa \|f_\alpha - f_\eta^*\|_\infty^\kappa$. This means that our assumption is less restrictive than requiring that the lower bound be tight.

Let moreover $\rho > 0$ be a parameter characterizing the complexity of the underlying set of boundary fragments and let $c_0$ be some constant. Denote by $\mathcal{F}_\rho$ a class of functions $f:[0,1]^d \to [0,1]$ satisfying the following condition: for every $f \in \mathcal{F}_\rho$ and every integer $m \leq L$ one has

$$(5.2) \qquad \min_{\alpha \,:\, m(\alpha) \leq m} \|f_\alpha - f\|_\infty \leq c_0 N_m^{-1/\rho}.$$

This is true for various smoothness classes (Sobolev, Hölder and certain Besov classes) with $1/\rho = \gamma/d$, where $\gamma$ is the regularity of the boundary $f$ (e.g., the number of bounded derivatives of $f$), and various bases $\{\psi_k\}$ [cf., e.g., Härdle, Kerkyacharian, Picard and Tsybakov (1998), Corollary 8.2 and Theorem 9.6].

Denote by $\mathcal{G}_\rho$ a class of boundary fragments $G = \{(s,t) \in \mathcal{X} : f(s) \geq t\}$ such that $f \in \mathcal{F}_\rho$.

THEOREM 2. *Suppose that Assumptions B and C are met. Then*

$$(5.3) \qquad \sup_{\eta \in \mathcal{H}_\kappa \,:\, G_\eta^* \in \mathcal{G}_\rho} [\mathbf{E}_\eta(R(\hat{G}_n)) - R(G_\eta^*)] = O\bigg(\bigg(\frac{\log^4 n}{n}\bigg)^{\kappa/(2\kappa+\rho-1)}\bigg),$$

*as $n \to \infty$.*

REMARK 2. For Hölder classes $\mathcal{F}_\rho$, the result of Theorem 2 is optimal up to a logarithmic factor [cf. Mammen and Tsybakov (1999) and Tsybakov (2004)]. Note that we cover here all values $\rho > 0$, thus extending the adaptive result of Tsybakov (2004) to $\rho \geq 1$ (i.e., to very irregular classes of boundaries). The case $\rho = 0$ can be also introduced: it corresponds to the assumption that (5.2) holds with 0 in the right-hand side. The class of functions $f$ thus defined is a Vapnik–Chervonenkis class, and it is easy to see that the rate in Theorem 2 in this case becomes $(n^{-1} \log^4 n)^{\kappa/(2\kappa-1)}$.

**6. Edge estimation.** In this section we consider the problem of estimation of the edge function $f_\eta^*$ such that $G_\eta^* = \{(s,t) \in \mathcal{X} : f_\eta^*(s) \geq t\}$, using the sample $\{(X_i, Y_i)\}_{i=1}^n$. The risk for this problem is defined by $\mathbf{E}(\mu_{d+1}(\hat{G}_n \triangle G_\eta^*)) = \mathbf{E}_\eta \|\hat{f}_n - f_\eta^*\|_1$ where $\hat{f}_n = f_{\hat{\alpha}_n}$ is the estimator of $f_\eta^*$ obtained by our method. Using the definition of $\mathcal{H}_\kappa$ we immediately get the following corollary of Theorem 2.

COROLLARY 1. *Suppose that Assumptions* B *and* C *are met. Then*

$$\sup_{\eta \in \mathcal{H}_\kappa : G_\eta^* \in \mathcal{G}_\rho} \mathbf{E}_\eta \|\hat{f}_n - f_\eta^*\|_1 = O\left(\left(\frac{\log^4 n}{n}\right)^{1/(2\kappa+\rho-1)}\right), \tag{6.1}$$

*as* $n \to \infty$.

Note that the setup of Corollary 1 is somewhat different from the standard problem of edge estimation as defined by Korostelev and Tsybakov (1993). In fact, it is in a sense more general because here the $X_i$'s are not supposed to be uniformly distributed on $[0,1]^d$ and the joint distribution of $(X, Y)$ is not supposed to follow a specified regression scheme. Also, the margin behavior is accounted for by the parameter $\kappa$. On the other hand, Corollary 1 deals only with binary images, $Y_i \in \{0, 1\}$, while Korostelev and Tsybakov (1993) allow $Y_i \in \mathbf{R}$, for instance, the model

$$Y_i = \mathbb{1}_{G^0}(X_i) + \xi_i, \qquad i = 1, \ldots, n, \tag{6.2}$$

where $\xi_i$ is a zero-mean random variable independent of $X_i$, and consider the problem of estimation of the edge function $f^0$ assuming that $G^0$ is a boundary fragment, $G^0 = \{(s,t) \in \mathcal{X} : f^0(s) \geq t\}$.

An important example covered by Corollary 1 is the model

$$Y_i = (1 + (2\mathbb{1}_{G^0}(X_i) - 1)\xi_i)/2, \qquad i = 1, \ldots, n, \tag{6.3}$$

where $\xi_i$ is a random variable independent of $X_i$ and taking values $-1$ and $1$ with probabilities $1 - p$ and $p$, respectively, $1/2 < p < 1$. In this model the observations $Y_i$ take values in $\{0, 1\}$ and they differ from the original (nonnoisy) image values $Y_i' = \mathbb{1}_{G^0}(X_i)$ because some values $Y_i'$ are switched from 0 to 1 and vice versa with probabilities $1 - p$ and $p$. This occurs, for example, if the image is transmitted through a binary channel. The aim is to estimate the edge function $f^0$ of the set $G^0$ assuming that $G^0$ is a boundary fragment.

It is easy to see that the regression function $\eta$ for the model (6.3) equals $\eta(x) = p\mathbb{1}_{G^0}(x) + (1-p)(1 - \mathbb{1}_{G^0}(x))$, which implies that the set $G^0$ is identical to $G_\eta^*$, and thus $f^0 = f_\eta^*$. Also, it is not hard to check that if the distribution of $X_i$'s is uniform on $[0,1]^{d+1}$ we have that $\eta \in \mathcal{H}_1$, and Corollary 1 applies with $\kappa = 1$.

Inspection of the proofs below shows that an analog of Corollary 1 also holds for the model (6.2) if one assumes that the random variables $Y_i$ are uniformly bounded. In this case only the constants in Lemma 4 and in the definition of $\lambda_n$ should be changed and the set $G^*$ should be indexed by the corresponding edge function $f$ rather than by the regression $\eta$, other elements of the construction remaining intact. This extension is quite obvious, and we do not pursue it here in more detail.

For $\kappa = 1$, Corollary 1 gives the rate $n^{-1/(\rho+1)}$, up to a logarithmic factor. As shown by Korostelev and Tsybakov (1993), this rate is optimal in a minimax sense when $\mathcal{F}_\rho$ is a Hölder class of functions and the model is (6.2) or (6.3). Barron, Birgé and Massart (1999) constructed adaptive estimators of the edge function in the model (6.2) with $d = 1$, $\kappa = 1$, $\rho \geq \rho_0 > 0$ using a penalization with a penalty that depends on the lower bound $\rho_0$ on $\rho$. They proved that for this particular case the optimal rate $n^{-1/(\rho+1)}$ is attained by their procedure. Corollary 1 extends these results, showing that our method allows adaptation to any complexity $\rho > 0$ in any dimension $d \geq 1$ and also adaptation to the margin $\kappa \geq 1$ which is necessary when we are not sure that the boundary is sharp, that is, when the regression function $\eta$ does not necessarily have a jump at the boundary. Assumption A or (5.1) gives a convenient characterization of nonsharpness of the boundary, and our penalized procedure allows us to adapt to the degree of non-sharpness.

**7. Proofs.** Before going into the technical details, let us first briefly explain our choice of class of sets as boundary fragments, and the choice of the penalty. When using boundary fragments, it is clear from (3.2) that the approximation of sets boils down to approximation of functions in $L_1$. We then use linear expansions, and need to relate the coefficients in these expansions to the penalty. This is done in Lemmas 1 and 2. Lemma 1 bounds the $L_1$-norm by $I(\cdot)$. Lemma 2 bounds $I(\cdot)$ by the $L_1$-norm when the number of levels is limited by $m$. The (block-wise) $\ell_{1/2}$ penalty ensures some important cancellations in the proof of Theorem 1. Its specific structure is less important in Lemmas 3 and 4, with Lemma 4 being a rather standard application of empirical process theory. Lemma 3 provides an upper bound for the *entropy with bracketing* (see the definition preceding Lemma 3) of the class of sets $G_\alpha \Delta G_{\alpha^*}$ with $\alpha$ varying, $\alpha^*$ fixed, and $I(\alpha - \alpha^*) \leq M$, $M > 0$. Lemma 4 is the consequence of the entropy result of Lemma 3 for the empirical process.

LEMMA 1. *Under Assumption* B *we have, for all* $\alpha \in \mathbf{R}^n$,
$$\|f_\alpha\|_1 \leq c_\psi I(\alpha). \tag{7.1}$$

PROOF. Using (3.4) we obtain
$$\|f_\alpha\|_1 = \left\| \sum_{j,l} \alpha_{j,l} \psi_{j,l} \right\|_1$$
$$\leq \sum_{j,l} |\alpha_{j,l}| \|\psi_{j,l}\|_1 \leq c_\psi \sum_{j,l} |\alpha_{j,l}| 2^{-dl/2}$$
$$= c_\psi \sum_{l=1}^{L} 2^{-dl} 2^{dl/2} \sum_j |\alpha_{j,l}|.$$

But clearly, for all $l$,

$$2^{dl/2} \sum_j |\alpha_{j,l}| = \left(2^{dl/4} \sqrt{\sum_j |\alpha_{j,l}|}\right)^2 \leq I(\alpha).$$

Hence,

$$\|f_\alpha\|_1 \leq c_\psi \sum_{l=1}^{L} 2^{-dl} I(\alpha) \leq c_\psi I(\alpha). \qquad \square$$

LEMMA 2. *Let $\alpha \in \mathbf{R}^n$ and let $N(\alpha)$ be defined in* (4.2). *Then under Assumption* B

(7.2) $$I(\alpha) \leq c_d c_\psi^2 N(\alpha) \|f_\alpha\|_1.$$

PROOF. The coefficient $\alpha_{j,l}$ is the inner product

$$\alpha_{j,l} = \int f_\alpha \psi_{j,l} \, d\mu_d,$$

so by (3.5),

$$\sum_{j \in I_l} |\alpha_{j,l}| \leq \int |f_\alpha| \sum_{j \in I_l} |\psi_{j,l}| \, d\mu_d$$

$$\leq c_\psi 2^{dl/2} \|f_\alpha\|_1.$$

This implies that for $m = m(\alpha)$, with $m(\alpha)$ given in (4.1),

$$\sqrt{I(\alpha)} = \sum_{l=1}^{m} 2^{dl/4} \sqrt{\sum_{j \in I_l} |\alpha_{j,l}|}$$

$$\leq \sum_{l=1}^{m} 2^{dl/2} \sqrt{c_\psi \|f_\alpha\|_1}$$

$$\leq \frac{2^{(m+1)d/2}}{2^{d/2} - 1} \sqrt{c_\psi \|f_\alpha\|_1}.$$

Next, by (3.6) and the definition (4.2) of $N(\alpha)$,

$$N(\alpha) = \sum_{l=1}^{m} |I_l| \geq c_\psi^{-1} \sum_{l=1}^{m} 2^{dl} \geq \frac{2^{(m+1)d}}{2c_\psi(2^d - 1)}.$$

Combining these inequalities we get the result. $\square$

DEFINITION 1. Let $\mathcal{Z} \subset L_p(\mathcal{S}, \nu)$ be a collection of functions on some measurable space $(\mathcal{S}, \nu)$, $1 \leq p \leq \infty$. For each $\delta > 0$, the $\delta$-covering number with bracketing $N_{B,p}(\delta, \mathcal{Z}, \nu)$ of $\mathcal{Z}$ is the smallest value of $N$ such that there exists a collection of pairs of functions $\{[z_j^L, z_j^U]_{j=1}^N\}$ that satisfies:

- $z_j^L \leq z_j^U$ and $\|z_j^U - z_j^L\|_p \leq \delta$ for all $j \in \{1, \ldots, N\}$ [with $\|\cdot\|_p$ being the $L_p(\mathcal{S}, \nu)$-norm],
- for each $z \in \mathcal{Z}$ there is a $j \in \{1, \ldots, N\}$ such that $z_j^L \leq z \leq z_j^U$.

The $\delta$-entropy with bracketing of $\mathcal{Z}$ is $H_{B,p}(\delta, \mathcal{Z}, \nu) = \log N_{B,p}(\delta, \mathcal{Z}, \nu)$.

DEFINITION 2. Let $\mathcal{Z}$ be a collection of bounded functions on $\mathcal{S}$. The $\delta$-covering number for the sup-norm, $N_\infty(\delta, \mathcal{Z})$, is the smallest number $N$ such that there are functions $\{z_j\}_{j=1}^N$ with for each $z \in \mathcal{Z}$,

$$\min_{j=1,\ldots,N} \sup_{s \in \mathcal{S}} |z(s) - z_j(s)| \leq \delta.$$

The $\delta$-entropy for the sup-norm is $H_\infty(\delta, \mathcal{Z}) = \log N_\infty(\delta, \mathcal{Z})$.

Note that when $\nu$ is a probability measure [cf. van de Geer (2000), page 17],

(7.3) $$H_{B,p}(\delta, \mathcal{Z}, \nu) \leq H_\infty(\delta/2, \mathcal{Z}), \delta > 0.$$

For a class $\mathcal{G}$ of subsets of $(\mathcal{X}, Q)$, we write $H_B(\delta, \mathcal{G}, Q) = H_{B,1}(\delta, \{\mathbb{1}_G : G \in \mathcal{G}\}, Q)$.

LEMMA 3. *Let $\alpha^* \in \mathbf{R}^n$ be fixed. For $0 < M \leq n$ define $\mathcal{G}^M = \{G = G_\alpha \triangle G_{\alpha^*} : \alpha \in \mathbf{R}^n, I(\alpha - \alpha^*) \leq M\}$. Suppose that Assumptions B and C are met. Then*

(7.4) $$H_B(\delta, \mathcal{G}^M, Q) \leq \frac{M}{\delta}\left(\frac{8q_0 c_\psi^2 \log n}{d}\right) \log\left(\frac{8q_0 c_\psi^2 n}{\delta d}\right),$$

*for all $0 < \delta \leq 1$.*

PROOF. Define $\mathcal{F}^M = \{f_\alpha : \alpha \in \mathbf{R}^n, I(\alpha) \leq M\}$. In view of Assumption C,

(7.5) $$H_B(q_0 \delta, \mathcal{G}^M, Q) \leq H_{B,1}(\delta, \mathcal{F}^M, \mu_d), \qquad \delta > 0.$$

This and (7.3) show that it is sufficient to bound $H_\infty(\cdot, \mathcal{F}^M)$.

Fix some $\delta > 0$. Our aim is now to bound the quantity $H_\infty((c_\psi^2 d^{-1} \log n)\delta, \mathcal{F}^M)$. To do this, note that one can construct a $(c_\psi^2 d^{-1} \log n)\delta$-net on $\mathcal{F}^M$ for the sup-norm in the following way. The elements of the net are $f_{\alpha'}$ where $\alpha'_{j,l}$ takes discretized values with step $\delta 2^{-dl/2}$. For every $\alpha_{j,l}$ define $\alpha'_{j,l}$ as the element closest to $\alpha_{j,l}$, of the $\delta 2^{-dl/2}$-net on the interval

$$[-M 2^{-dl/2}, M 2^{-dl/2}].$$

Note that this interval contains all admissible values of $\alpha_{j,l}$ since $|\alpha_{j,l}| \leq M 2^{-dl/2}$, $\forall j, l$ for all $\alpha$ such that $I(\alpha) \leq M$. With this definition of $\alpha'_{j,l}$ we have $|\alpha_{j,l} - \alpha'_{j,l}| \leq \delta 2^{-dl/2}$, and thus

$$\sup_{s \in [0,1]^d} |f_\alpha(s) - f_{\alpha'}(s)|$$

$$\leq \sum_{l=1}^{L} \sup_{s \in [0,1]^d} \sum_{j \in I_l} |\alpha_{j,l} - \alpha'_{j,l}| |\psi_{j,l}(s)|$$

$$\leq \delta \sum_{l=1}^{L} 2^{-dl/2} \sup_{s \in [0,1]^d} \sum_{j \in I_l} |\psi_{j,l}(s)| \leq L c_\psi \delta \leq (c_\psi^2 d^{-1} \log n) \delta,$$

where we have used Assumption B for the last two inequalities. Thus we have proved that the above construction gives in fact a $(c_\psi^2 d^{-1} \log n)\delta$-net on $\mathcal{F}^M$ for the sup-norm.

Let us now evaluate the cardinality of this net. This will be based on the following three observations.

OBSERVATION 1. For every $\alpha$ such that $I(\alpha) \leq M$ there exist at most $M/\delta$ indices $k = (j, l)$ such that $|\alpha_{j,l}| > \delta 2^{-dl/2}$. To show this, define

$$N_l(\alpha) = |\{j \in I_l : |\alpha_{j,l}| > \delta 2^{-dl/2}\}|, \qquad l = 1, \ldots, L.$$

Then

$$\sqrt{M} \geq \sqrt{I(\alpha)} \geq \sum_{l=1}^{L} 2^{dl/4} \sqrt{\sum_{|\alpha_{j,l}| > \delta 2^{-dl/2}} |\alpha_{j,l}|} \geq \sqrt{\delta} \sum_{l=1}^{L} \sqrt{N_l(\alpha)}.$$

Hence

$$\sum_{l=1}^{L} \sqrt{N_l(\alpha)} \leq \sqrt{\frac{M}{\delta}},$$

and so

$$\sum_{l=1}^{L} N_l(\alpha) \leq \frac{M}{\delta}.$$

OBSERVATION 2. For each $j$ and $l$, we can approximate the interval $\{|\alpha_{j,l}| \leq M 2^{-dl/2}\}$ by a set of cardinality at most

$$\frac{2M}{\delta} + 1$$

such that each coefficient $\alpha_{j,l}$ is approximated to within the distance $\delta 2^{-dl/2}$.

OBSERVATION 3. The number of different ways to choose $\leq M/\delta$ nonzero coefficients out of $n$ is

$$\sum_{0 \leq N \leq \min\{M/\delta, n\}} \binom{n}{N} \leq (n+1)^{M/\delta}$$

[see, e.g., Devroye, Györfi and Lugosi (1996), page 218].

It follows from Observation 3 that there exist at most $(n+1)^{M/\delta}$ possibilities to choose the sets of nonzero coordinates of the vectors $\alpha'$ belonging to the net. For each of these possibilities the discretization is performed on each of the nonzero coordinates, which gives at most

$$\left(\frac{2M}{\delta} + 1\right)^{M/\delta}$$

new possibilities in view of Observations 1 and 2. Thus, the cardinality of the considered $(c_\psi^2 d^{-1} \log n)\delta$-net on $\mathcal{F}^M$ is bounded by

$$(n+1)^{M/\delta}\left(\frac{2M}{\delta} + 1\right)^{M/\delta},$$

which implies

(7.6) $\quad H_\infty((c_\psi^2 d^{-1} \log n)\delta, \mathcal{F}^M) \leq \dfrac{M}{\delta}\left(\log\left(\dfrac{2M}{\delta} + 1\right) + \log(n+1)\right).$

In view of (7.5) this yields

$$H_B(\delta, \mathcal{G}^M, Q) \leq \frac{M}{\delta}\left(\frac{2q_0 c_\psi^2 \log n}{d}\right)\left[\log\left(\frac{4q_0 c_\psi^2 M \log n}{\delta d} + 1\right) + \log(n+1)\right]$$

$$\leq \frac{M}{\delta}\left(\frac{2q_0 c_\psi^2 \log n}{d}\right)\left[\log\left(\frac{4q_0 c_\psi^2 n^2}{\delta d} + 1\right) + \log(n+1)\right]$$

since $M \log n \leq n \log n \leq n^2$. Continuing with this bound, we arrive at

$$H_B(\delta, \mathcal{G}^M, Q) \leq \frac{M}{\delta}\left(\frac{4q_0 c_\psi^2 \log n}{d}\right)\log\left(\frac{4q_0 c_\psi^2 n^2}{\delta d} + 1\right)$$

$$\leq \frac{M}{\delta}\left(\frac{4q_0 c_\psi^2 \log n}{d}\right)\log\left(\frac{8q_0 c_\psi^2 n^2}{\delta d}\right)$$

$$\leq \frac{M}{\delta}\left(\frac{8q_0 c_\psi^2 \log n}{d}\right)\log\left(\frac{8q_0 c_\psi^2 n}{\delta d}\right). \qquad \square$$

Now we turn to the empirical process

(7.7) $\qquad \nu_n(\alpha) = \sqrt{n}(R_n(G_\alpha) - R(G_\alpha)), \qquad \alpha \in \mathbf{R}^n.$

LEMMA 4. *Let Assumptions B and C hold. Then there exists a universal constant $C$ such that for $n \geq 8q_0 c_\psi^2$ we have, for all $\alpha^* \in \mathbf{R}^n$,*

(7.8)
$$\mathbf{P}\left(\sup_{\alpha \in \mathbf{R}^n} \frac{|\nu_n(\alpha) - \nu_n(\alpha^*)|}{\sqrt{I(\alpha - \alpha^*)} + \sqrt{\log^4 n/n}} > C\sqrt{\frac{q_0 c_\psi^2 \log^4 n}{d}}\right)$$
$$\leq C \exp\left[-\frac{c_\psi \log^4 n}{C^2 d}\right].$$

PROOF. We will apply Theorem 5.11 in van de Geer (2000) which, translated to our situation, says the following. Let
$$h_\alpha(X, Y) = (Y - \mathbb{1}_{G_\alpha}(X))^2 - (Y - \mathbb{1}_{G_{\alpha^*}}(X))^2$$
and
$$\mathcal{H}^M = \{h_\alpha : I(\alpha - \alpha^*) \leq M\}.$$

Also, let $R^2 \leq 1$ satisfy
$$\sup_{h \in \mathcal{H}^M} \int h^2 \, dP \leq R^2,$$
where $P$ is the law of $(X, Y)$. Then Theorem 5.11 in van de Geer (2000) gives that for some universal constant $C_0$, and for all $a$ satisfying both $a \leq \sqrt{n}R^2$ and
$$a \geq C_0 \left(\int_{a/(C_0\sqrt{n})}^1 H_{B,2}^{1/2}(u, \mathcal{H}^M, P) \, du \vee R\right)$$
one has

(7.9) $\mathbf{P}\left(\sup_{\alpha \in \mathbf{R}^n : I(\alpha - \alpha^*) \leq M} |\nu_n(\alpha) - \nu_n(\alpha^*)| > a\right) \leq C_0 \exp\left[-\frac{a^2}{C_0^2 R^2}\right].$

To apply this result, note first that

(7.10) $|(Y - \mathbb{1}_{G_\alpha}(X))^2 - (Y - \mathbb{1}_{G_{\alpha'}}(X))^2| = |\mathbb{1}_{G_\alpha}(X) - \mathbb{1}_{G_{\alpha'}}(X)|.$

We therefore get
$$\sup_{h \in \mathcal{H}^M} \int h^2 \, dP = \sup_{G \in \mathcal{G}^M} Q(G),$$
where $\mathcal{G}^M$ be defined as in Lemma 3. Hence by Lemma 1, Assumption C and (3.2) we may take
$$R^2 = q_0 c_\psi M \wedge 1.$$

Moreover, again by (7.10),
$$H_{B,2}(\delta, \mathcal{H}^M, P) = H_{B,2}(\delta, \{\mathbb{1}_G : G \in \mathcal{G}^M\}, Q) = H_B(\delta^2, \mathcal{G}^M, Q), \qquad \delta > 0.$$

Using Lemma 3, for any $a \leq \sqrt{n}R^2$, $\log^4 n/n \leq M \leq n$ and $n \geq 8q_0 c_\psi^2$, we get the bound

$$\int_{a/\sqrt{n}}^1 H_{B,2}^{1/2}(u, \mathcal{H}^M, P) \, du \leq c' \sqrt{\frac{q_0 c_\psi^2 M \log n}{d}} \log^3\left(\frac{n^{5/2}}{a^2}\right),$$

where $c'$ is a universal constant. We therefore can take

$$a = c\sqrt{\frac{q_0 c_\psi^2 M \log^4 n}{d}},$$

with an appropriate universal constant $c$. Insert this value for $a$ and the value of $R$ in (7.9) to find that for $\log^4 n/n \leq M \leq n$, and trivially also for $M > n$,

(7.11)
$$\mathbf{P}\left(\sup_{\alpha \in \mathbf{R}^n : I(\alpha - \alpha^*) \leq M} |\nu_n(\alpha) - \nu_n(\alpha^*)| > c\sqrt{\frac{q_0 c_\psi^2 M \log^4 n}{d}}\right)$$
$$\leq C_0 \exp\left[-\frac{c^2 c_\psi \log^4 n}{C_0^2 d}(M \vee 1)\right].$$

The result now follows from the peeling device as, for example, explained in Section 5.3 of van de Geer (2000). The argument is then as follows. We have

$$\mathbf{P}\left(\sup_{\alpha \in \mathbf{R}^n} \frac{|\nu_n(\alpha) - \nu_n(\alpha^*)|}{\sqrt{I(\alpha - \alpha^*)} + \sqrt{\log^4 n/n}} > C\sqrt{\frac{q_0 c_\psi^2 \log^4 n}{d}}\right)$$
$$\leq \mathbf{P}\left(\sup_{I(\alpha - \alpha^*) \leq 1} \frac{|\nu_n(\alpha) - \nu_n(\alpha^*)|}{\sqrt{I(\alpha - \alpha^*)} + \sqrt{\log^4 n/n}} > C\sqrt{\frac{q_0 c_\psi^2 \log^4 n}{d}}\right)$$
$$+ \mathbf{P}\left(\sup_{I(\alpha - \alpha^*) > 1} \frac{|\nu_n(\alpha) - \nu_n(\alpha^*)|}{\sqrt{I(\alpha - \alpha^*)} + \sqrt{\log^4 n/n}} > C\sqrt{\frac{q_0 c_\psi^2 \log^4 n}{d}}\right)$$
$$= \mathbf{P}_I + \mathbf{P}_{II}.$$

Furthermore, for $j_0$ the integer such that $2^{-j_0} \leq \log^4 n/n < 2^{-j_0+1}$, we find

$$\mathbf{P}_I \leq \sum_{j=0}^{j_0} \mathbf{P}\left(\sup_{I(\alpha - \alpha^*) \leq 2^{-j}} |\nu_n(\alpha) - \nu_n(\alpha^*)| > \frac{C}{2}\sqrt{\frac{q_0 c_\psi^2 2^{-j} \log^4 n}{d}}\right) = \sum_{j=0}^{j_0} \mathbf{P}_{I,j}.$$

Similarly,

$$\mathbf{P}_{II} \leq \sum_{j=1}^{\infty} \mathbf{P}\left(\sup_{I(\alpha - \alpha^*) \leq 2^j} |\nu_n(\alpha) - \nu_n(\alpha^*)| > \frac{C}{2}\sqrt{\frac{q_0 c_\psi^2 2^j \log^4 n}{d}}\right) = \sum_{j=1}^{\infty} \mathbf{P}_{II,j}.$$

The theorem then follows by choosing $C$ appropriately and applying (7.11) to each of the $\mathbf{P}_{I,j}$, $j = 0, \ldots, j_0$, and $\mathbf{P}_{II,j}$, $j = 1, 2, \ldots$. □

LEMMA 5. *For any positive $v, t$ and any $\kappa \geq 1$, $\delta > 0$ we have*
$$vt^{1/(2\kappa)} \leq (\delta/2)t + c_\kappa \delta^{-1/(2\kappa-1)} v^{2\kappa/(2\kappa-1)}$$
*where $c_\kappa = (2\kappa - 1)/(2\kappa)\kappa^{-1/(2\kappa-1)}$.*

PROOF. By the concavity of the log-function, we have for positive $a$, $b$, $x$ and $y$, with $1/x + 1/y = 1$,
$$\log(ab) = \frac{1}{x}\log(a^x) + \frac{1}{y}\log(b^y) \leq \log\left(\frac{1}{x}a^x + \frac{1}{y}b^y\right)$$

or

$$ab \leq \frac{1}{x}a^x + \frac{1}{y}b^y.$$

The lemma is obtained when we choose
$$a = v(\kappa\delta)^{-1/(2\kappa)}, \qquad b = (\kappa\delta t)^{1/(2\kappa)}, \qquad x = \frac{2\kappa}{2\kappa - 1}, \qquad y = 2\kappa. \qquad □$$

We now come to the proof of the main theorem. This proof follows the lines of Loubes and van de Geer (2002) [see also van de Geer (2003)].

PROOF OF THEOREM 1. Fix an arbitrary $\alpha^* \in \mathbf{R}^n$. (We stress here that $\alpha^*$ is just a notation and need not be related in any sense to the Bayes rule $G^*$.) Let $\Xi$ be the random event

$$(7.12) \quad \Xi = \left\{|\nu_n(\hat{\alpha}_n) - \nu_n(\alpha^*)|/\sqrt{n} \leq \lambda_n\sqrt{I(\hat{\alpha}_n - \alpha^*)} + \lambda_n\sqrt{\frac{\log^4 n}{n}}\right\}.$$

By Lemma 4, for $n$ sufficiently large,
$$\mathbf{P}(\Xi) \geq 1 - C\exp\left[-\frac{c_\psi \log^4 n}{C^2 d}\right].$$

So we only need to consider what happens on the set $\Xi$. The definition of $\hat{\alpha}_n$ implies
$$R_n(G_{\hat{\alpha}_n}) + \lambda_n\sqrt{I(\hat{\alpha}_n)} \leq R_n(G_{\alpha^*}) + \lambda_n\sqrt{I(\alpha^*)},$$

which may be rewritten in the form

$$(7.13) \quad R(G_{\hat{\alpha}_n}) \leq -[\nu_n(\hat{\alpha}_n) - \nu_n(\alpha^*)]/\sqrt{n} - \lambda_n[\sqrt{I(\hat{\alpha}_n)} - \sqrt{I(\alpha^*)}] + R(G_{\alpha^*}).$$

Hence on $\Xi$ we get

$$R(G_{\hat{\alpha}_n}) \le \lambda_n\sqrt{I(\hat{\alpha}_n - \alpha^*)} - \lambda_n[\sqrt{I(\hat{\alpha}_n)} - \sqrt{I(\alpha^*)}\,] + R(G_{\alpha^*}) + \lambda_n\sqrt{\frac{\log^4 n}{n}}.$$

Let $m^* = m(\alpha^*)$, and let, for any $\alpha$,

$$\sqrt{I^{(1)}(\alpha)} = \sum_{l=1}^{m^*} 2^{dl/4}\sqrt{\sum_{j \in I_l}|\alpha_{j,l}|},$$

$$\sqrt{I^{(2)}(\alpha)} = \sum_{l=m^*+1}^{L} 2^{dl/4}\sqrt{\sum_{j \in I_l}|\alpha_{j,l}|}.$$

Since $I^{(2)}(\alpha - \alpha^*) = I^{(2)}(\alpha)$, we now find

$$R(G_{\hat{\alpha}_n}) \le \lambda_n\sqrt{I^{(1)}(\hat{\alpha}_n - \alpha^*)} + \lambda_n\sqrt{I^{(2)}(\hat{\alpha}_n)} - \lambda_n[\sqrt{I^{(1)}(\hat{\alpha}_n)} - \sqrt{I^{(1)}(\alpha^*)}\,]$$

$$- \lambda_n\sqrt{I^{(2)}(\hat{\alpha}_n)} + R(G_{\alpha^*}) + \lambda_n\sqrt{\frac{\log^4 n}{n}}$$

$$= \lambda_n\sqrt{I^{(1)}(\hat{\alpha}_n - \alpha^*)} - \lambda_n[\sqrt{I^{(1)}(\hat{\alpha}_n)} - \sqrt{I^{(1)}(\alpha^*)}\,]$$

$$+ R(G_{\alpha^*}) + \lambda_n\sqrt{\frac{\log^4 n}{n}}.$$

Since for any $a, b \in \mathbf{R}$, $\sqrt{|a|} - \sqrt{|b|} \le \sqrt{|a-b|}$, we arrive at

(7.14) $\quad R(G_{\hat{\alpha}_n}) \le 2\lambda_n\sqrt{I^{(1)}(\hat{\alpha}_n - \alpha^*)} + R(G_{\alpha^*}) + \lambda_n\sqrt{\dfrac{\log^4 n}{n}}.$

Therefore, using a straightforward modification of Lemma 2 (basically replacing there $I$ by $I^{(1)}$), we obtain

$$R(G_{\hat{\alpha}_n}) \le 2\lambda_n\sqrt{c_d c_\psi^2 N^* \|f_{\hat{\alpha}_n - \alpha^*}\|_1} + R(G_{\alpha^*}) + \lambda_n\sqrt{\frac{\log^4 n}{n}},$$

where $N^* = N(\alpha^*) = \sum_{l=1}^{m^*}|I_l|$. By Assumption C and (3.2),

$$\|f_{\hat{\alpha}_n - \alpha^*}\|_1 \le q_0 Q(G_{\hat{\alpha}_n} \triangle G_{\alpha^*}).$$

We therefore get

$$R(G_{\hat{\alpha}_n}) \le 2\lambda_n\sqrt{c_d q_0 c_\psi^2 N^* Q(G_{\hat{\alpha}_n} \triangle G_{\alpha^*})} + R(G_{\alpha^*}) + \lambda_n\sqrt{\frac{\log^4 n}{n}}.$$

Subtracting $R(G^*)$ from both sides of this inequality, and denoting $d(G, G^*) = R(G) - R(G^*)$, we obtain

$$d(G_{\hat{\alpha}_n}, G^*) \leq 2\lambda_n \sqrt{c_d q_0 c_\psi^2 N^* Q(G_{\hat{\alpha}_n} \triangle G_{\alpha^*})} \tag{7.15}$$

$$+ d(G_{\alpha^*}, G^*) + \lambda_n \sqrt{\frac{\log^4 n}{n}}.$$

But then, by the triangle inequality and $\sqrt{a+b} \leq \sqrt{a} + \sqrt{b}$, $a, b \geq 0$, we get

$$d(G_{\hat{\alpha}_n}, G^*) \leq 2\lambda_n \sqrt{c_d q_0 c_\psi^2 N^*} [\sqrt{Q(G_{\hat{\alpha}_n} \triangle G^*)} + \sqrt{Q(G_{\alpha^*} \triangle G^*)}]$$

$$+ d(G_{\alpha^*}, G^*) + \lambda_n \sqrt{\frac{\log^4 n}{n}}$$

$$\leq 2\lambda_n \sqrt{c_d q_0 c_\psi^2 \sigma_0^{1/\kappa} N^*} [d^{1/(2\kappa)}(G_{\hat{\alpha}_n}, G^*) + d^{1/(2\kappa)}(G_{\alpha^*}, G^*)]$$

$$+ d(G_{\alpha^*}, G^*) + \lambda_n \sqrt{\frac{\log^4 n}{n}},$$

where in the last inequality we invoked Assumption A. Now we apply Lemma 5 with, respectively, $t = d(G_{\hat{\alpha}_n}, G^*)$ and $t = d(G_{\alpha^*}, G^*)$, to get

$$d(G_{\hat{\alpha}_n}, G^*) \leq (\delta/2)[d(G_{\hat{\alpha}_n}, G^*) + d(G_{\alpha^*}, G^*)]$$

$$+ 2c_\kappa \delta^{-1/(2\kappa-1)} (4 c_d q_0 c_\psi^2 \sigma_0^{1/\kappa} \lambda_n^2 N^*)^{\kappa/(2\kappa-1)}$$

$$+ d(G_{\alpha^*}, G^*) + \lambda_n \sqrt{\frac{\log^4 n}{n}},$$

which, together with the inequalities $(1 + \delta/2)/(1 - \delta/2) \leq (1 + \delta)^2$ and $1/(1 - \delta/2) \leq 2$, which are valid for $\delta \in (0, 1]$, implies, that on the event $\Xi$ we have

$$R(G_{\hat{\alpha}_n}) - R(G^*)$$

$$\leq (1 + \delta)^2 \{R(G_{\alpha^*}) - R(G^*) + \delta^{-1/(2\kappa-1)} V_n(N(\alpha^*))\} + 2\lambda_n \sqrt{\frac{\log^4 n}{n}}.$$

Hence

$$\mathbf{P}\bigg(R(G_{\hat{\alpha}_n}) - R(G^*) > (1 + \delta)^2 \{R(G_{\alpha^*}) - R(G^*)$$

$$+ \delta^{-1/(2\kappa-1)} V_n(N(\alpha^*))\} + 2\lambda_n \sqrt{\frac{\log^4 n}{n}}\bigg)$$

$$\leq C \exp\left[-\frac{c_\psi \log^4 n}{C^2 d}\right].$$

Since $\alpha^*$ was chosen arbitrarily this holds in fact for all $\alpha^*$. Because a distribution function is right continuous, we now have shown that also

$$\mathbf{P}\bigg(R(G_{\hat{\alpha}_n}) - R(G^*) > (1+\delta)^2 \inf_{\alpha \in \mathbf{R}^n} \{R(G_\alpha) - R(G^*)$$

$$+ \delta^{-1/(2\kappa-1)} V_n(N(\alpha))\} + 2\lambda_n \sqrt{\frac{\log^4 n}{n}}\bigg)$$

$$\leq C \exp\left[-\frac{c_\psi \log^4 n}{C^2 d}\right]. \qquad \square$$

PROOF OF THEOREM 2. For $\eta \in \mathcal{H}_\kappa$, $G_\eta^* \in \mathcal{G}_\rho$, we have
$$R(G_\alpha) - R(G_\eta^*) + V_n(N(\alpha)) \leq \sigma_0 q_0^\kappa \|f_\alpha - f_\eta^*\|_\infty^\kappa + V_n(N(\alpha)),$$

so that

(7.16) $$\inf_{\alpha:\, m(\alpha) \leq m} \{R(G_\alpha) - R(G_\eta^*) + V_n(N(\alpha))\} \leq \sigma_0 q_0^\kappa c_0^\kappa N_m^{-\kappa/\rho} + V_n(N_m)$$
$$= z(N_m),$$

where
$$z(t) = \sigma_0 q_0^\kappa c_0^\kappa t^{-\kappa/\rho} + V_n(t), \qquad t > 0.$$

Now minimizing $z(t)$ over all $t > 0$ gives

$$t \asymp \left(\frac{n}{\log^4 n}\right)^{\rho/(2\kappa+\rho-1)} := \tilde{t},$$

since $V_n(N) \asymp (Nn^{-1} \log^4 n)^{\kappa/(2\kappa-1)}$. Let $\tilde{m}$ be the smallest integer such that

$$N_{\tilde{m}-1} \leq \tilde{t} \leq N_{\tilde{m}}.$$

It is not difficult to see, using (3.6) and (3.7), that

$$N_{\tilde{m}} - \tilde{t} \leq c_\psi^2 2^{2d} (\tilde{t} + 1).$$

Inserting $N_{\tilde{m}}$ in the right-hand side of (7.16) therefore gives

$$\inf_{\alpha:m(\alpha)\leq \tilde{m}} \{R(G_\alpha) - R(G_\eta^*) + V_n(N(\alpha))\} \leq z(N_{\tilde{m}}) \asymp \left(\frac{\log^4 n}{n}\right)^{\kappa/(2\kappa+\rho-1)}.$$

Note finally that the constants in Theorem 1 depend only on $d$, $\kappa$, $\sigma_0$, $q_0$ and $c_\psi$, so that the result of Theorem 2 follows easily. $\square$

REMARK. When this paper was finished we learned from Vladimir Koltchinskii that he found another penalized classifier that adaptively attains fast optimal rates [Koltchinskii (2003)]. His method is different from ours and uses randomization and local Rademacher complexities.

# REFERENCES

Audibert, J.-Y. (2004). Aggregated estimators and empirical complexity for least squares regression. *Ann. Inst. H. Poincaré Probab. Statist.* **40** 685–736. MR2096215

Barron, A., Birgé, L. and Massart, P. (1999). Risk bounds for model selection via penalization. *Probab. Theory Related Fields* **113** 301–413. MR1679028

Bartlett, P. L., Jordan, M. I. and McAuliffe, J. D. (2003). Convexity, classification and risk bounds. Technical Report 638, Dept. Statistics, Univ. California, Berkeley.

Blanchard, G., Lugosi, G. and Vayatis, N. (2003). On the rate of convergence of regularized boosting classifiers. *J. Mach. Learn. Res.* **4** 861–894. MR2076000

Cavalier, L. and Tsybakov, A. B. (2001). Penalized blockwise Stein's method, monotone oracles and sharp adaptive estimation. *Math. Methods Statist.* **10** 247–282. MR1867161

DeVore, R. A. and Lorentz, G. G. (1993). *Constructive Approximation.* Springer, Berlin. MR1261635

Devroye, L., Györfi, L. and Lugosi, G. (1996). *A Probabilistic Theory of Pattern Recognition.* Springer, New York. MR1383093

Härdle, W., Kerkyacharian, G., Picard, D. and Tsybakov, A. (1998). *Wavelets, Approximation and Statistical Applications. Lecture Notes in Statist.* **129**. Springer, New York. MR1618204

Koltchinskii, V. (2001). Rademacher penalties and structural risk minimization. *IEEE Trans. Inform. Theory* **47** 1902–1914. MR1842526

Koltchinskii, V. (2003). Local Rademacher complexities and oracle inequalities in risk minimization. Preprint.

Koltchinskii, V. and Panchenko, D. (2002). Empirical margin distributions and bounding the generalization error of combined classifiers. *Ann. Statist.* **30** 1–50. MR1892654

Korostelev, A. P. and Tsybakov, A. B. (1993). *Minimax Theory of Image Reconstruction. Lecture Notes in Statist.* **82**. Springer, New York. MR1226450

Loubes, J.-M. and van de Geer, S. (2002). Adaptive estimation with soft thresholding penalties. *Statist. Neerlandica* **56** 453–478. MR2027536

Lugosi, G. and Wegkamp, M. (2004). Complexity regularization via localized random penalties. *Ann. Statist.* **32** 1679–1697. MR2089138

Mammen, E. and Tsybakov, A. B. (1999). Smooth discrimination analysis. *Ann. Statist.* **27** 1808–1829. MR1765618

Schölkopf, B. and Smola, A. (2002). *Learning with Kernels*. MIT Press, Cambridge, MA.

Tsybakov, A. B. (2004). Optimal aggregation of classifiers in statistical learning. *Ann. Statist.* **32** 135–166. MR2051002

van de Geer, S. (2000). *Empirical Processes in M-Estimation*. Cambridge Univ. Press.

van de Geer, S. (2003). Adaptive quantile regression. In *Recent Advances and Trends in Nonparametric Statistics* (M. G. Akritas and D. N. Politis, eds.) 235–250. North-Holland, Amsterdam.

Vapnik, V. N. (1998). *Statistical Learning Theory*. Wiley, New York. MR1641250

Laboratoire de Probabilités
et Modèles Aléatoires
Université Paris VI
4 Place Jussieu
Case 188
F-75252 Paris Cédex 05
France
e-mail: tsybakov@ccr.jussieu.fr

Mathematical Institute
University of Leiden
P.O. Box 9512
2300 RA Leiden
The Netherlands
e-mail: geer@math.leidenuniv.nl